\documentclass[hidelinks,onefignum,onetabnum]{siamart220329}

\usepackage{graphicx}
\usepackage{epstopdf}
\ifpdf
  \DeclareGraphicsExtensions{.eps,.pdf,.png,.jpg}
\else
  \DeclareGraphicsExtensions{.eps}
\fi

\usepackage{amsmath,amssymb,amsfonts,amscd}
\usepackage[all,cmtip]{xy}
\usepackage{enumerate}
\usepackage{enumitem}
\usepackage{amssymb, latexsym, amsmath}
\usepackage{url}
\usepackage{nicematrix}
\usepackage{tikz}
\usepackage{caption, subcaption}
\usepackage{mathtools}
\usepackage[normalem]{ulem}
\usepackage{tikz}
\usepackage{pgfplots}

\usepackage[normalem]{ulem}
\newcommand{\GN}{generalized Nyström}
\newcommand{\Ny}{Nyström}
\newcommand{\GNa}{generalized Nyström approximation}

\newcommand{\F}{\mathbb{F}}
\newcommand{\jcol}[2]{#1_{:, #2}} 
\newcommand{\jrow}[2]{#1_{#2, :}} 
\newcommand{\jk}[3]{[#1]_{#2, #3}} 
\newcommand{\njcol}[2]{#1_{:, -#2}} 
\newcommand{\njrow}[2]{#1_{-#2, :}} 
\newcommand{\njk}[3]{#1_{-#2, -#3}} 
\newcommand{\hb}[2]{{h}_{#1}} 

\newcommand{\boxplotscale}{.6}

\definecolor{cSVD}{RGB}{0, 0, 0}
\definecolor{cGN}{RGB}{227, 178, 60}
\definecolor{cLTO}{RGB}{0, 20, 83}
\definecolor{cLPO}{RGB}{251, 97, 7}
\definecolor{cLCO}{RGB}{247, 37, 133}

\pgfplotsset{
    numformat single/.style={pgf/number format/#1},
    numformat/.style={numformat single/.list={#1}},
    every axis/.append style={
         width = .45 \textwidth,
        height = .4 \textwidth,
        every tick/.style = {xtick pos = left},
        title style={font=\bfseries},
        label style={font=\small},
        ticklabel style={font=\footnotesize},
        legend style={font=\tiny},
        scaled y ticks = false,
    },
    every axis plot/.append style={
    every mark/.append style={solid,  fill opacity = .7, draw opacity = .7, thick},
    gn_line/.style = {mark=triangle*, mark size = 3pt, color = cGN, dashed, line width = .8mm},
    lto_line/.style={mark=star, mark size = 2pt, color = cLTO, dashed, thick},
    blto_line/.style={mark=o, mark size = 2pt, color = cLTO, opacity = .5, thick},
    svds_line/.style={mark=square, mark size = .5pt, color = cSVD, draw opacity = 0.2, thin},
    lpo_line/.style={mark=star, mark size = 2pt, color = cLPO, dashed, thick},
    blpo_line/.style={mark=o, mark size = 2pt, color = cLPO, opacity = .5, thick},
    lco_line/.style={mark=star, mark size = 2pt, color = cLCO, dashed, thick},
    blco_line/.style={mark=o, mark size = 2pt, color = cLCO, opacity = .5, draw opacity = .5, thick},
    }
}

\newsiamremark{remark}{Remark}
\newsiamremark{hypothesis}{Hypothesis}
\crefname{hypothesis}{Hypothesis}{Hypotheses}
\newsiamthm{claim}{Claim}

\headers{{Leave-one-out estimator for \GN}}{L. Lazzarino, K. Pearce, N. Pritchard}

\title{Efficient error estimators for Generalized Nyström }
\author{Lorenzo Lazzarino\thanks{Mathematical Institute, University of Oxford, Oxford, OX2 6GG, UK  (\email{lorenzo.lazzarino@maths.ox.ac.uk}, \email{pritchard@maths.ox.ac.uk}).}, 
\and Katherine J. Pearce\thanks{Oden Institute, University of Texas at Austin, US (\email{katherine.pearce@austin.utexas.edu}).}
 \and Nathaniel Pritchard\footnotemark[1]}

\usepackage{amsopn}

\begin{document}

\maketitle

\frenchspacing

\begin{abstract}
Randomized algorithms in numerical linear algebra have proven to be effective in ameliorating issues of scalability when working with large matrices, efficiently producing accurate low-rank approximations. 
A key remaining challenge, however, is to efficiently assess the approximation accuracy of randomized methods without additional expensive matrix accesses. 
Recent work has addressed this issue by deriving fast leave-one-out error estimators for the randomized SVD and Nyström decomposition, enabling accurate error estimation with no additional matrix accesses. 
In this work, we extend the leave-one-out framework to the generalized Nyström decomposition, an approach that can be applied to general rectangular matrices. 
We do this by deriving three new leave-one-out error estimators and validating their effectiveness through numerical experiments.
\end{abstract}



\section{Introduction} \label{sec:Intro}

Low-rank approximation is an essential tool for accelerating computations in a myriad of applications across scientific computing. 
Classically, many low-rank approximations are computed as truncated singular value decompositions (SVDs) because the rank-$k$ truncated SVD yields the optimal  approximation error among all rank-$k$ approximations \cite{Eckart1936}.
Unfortunately, the asymptotic complexity of the SVD can make it an undesirable choice for large problems: given an $m \times m$ matrix $A$, the rank-$k$ truncated SVD has a complexity of $O(m^2k)$.

Recent advances in randomized numerical linear algebra have led to more efficient algorithms to compute low-rank decompositions, including the randomized SVD \cite{halkoFindingStructureRandomness2011}, interpolative \cite{dong2024robust, dong2023simpler, pearce2025adaptive}, CUR  \cite{mahoney2009cur, pritchard2025fast}, Cholesky \cite{chen2020efficient}, \Ny{} \cite{GittensMahoney2016}, and \GN{} \cite{nakatsukasaFastStableRandomized2020, tropp2017practical} decompositions.
The benefits of randomization are especially salient for large problems, with better asymptotic scaling and superior performance over classical algorithms, particularly in instances where a matrix can only be accessed through matrix-vector products.
However, once a low-rank approximation of a matrix has been computed, the problem remains of how to compute the approximation error. 
Na\"ively, this would require an additional access to the full matrix $A$, an often impractical or even impossible requirement in many scenarios, e.g. streaming data. Randomized methods have been proposed to estimate the error without additional full matrix accesses, (see \cite[Section 12.1]{martinssonRandomizedNumericalLinear2020a}) but still require additional matrix-vector products to form an independent random sample. 
However, to obtain better performance, it is crucial that the sample used to compute the low-rank approximation be independent of the sample used to estimate its error.
    
In \cite{epperly2024efficient}, a new error estimation approach is proposed that avoids these expensive full matrix accesses, called the leave-one-out estimator. 
In statistical settings, the leave-one-out estimator is frequently used to determine a model's generalization error, how well the model predicts on unseen data \cite{wasserman2010alla}.
The leave-one-out estimator provides this estimate of the error without requiring new data, by fitting a model on each possible partition of the data that includes all data points except for one, and computing the error at each excluded data point. After averaging the computed errors, the leave-one-out method provides an error estimate of the model if all data points were used for fitting. 
This approach naturally extends to randomized low-rank approximation methods; namely, in randomized low-rank approximation, a small set of test vectors (the data) are formed, or sketched, and used to compute an approximation (the model).  
Repeatedly computing these approximations over all data partitions is very expensive, so for this approach to be practical, a fast formula must be derived.

For the randomized SVD and \Ny{} decomposition, a fast leave-one-out error estimation formula is derived in \cite{epperly2024efficient}, enabling efficient estimation of the rank-$(k-1)$ approximation error using only the information obtained from a random sketch with $k$ columns.
When the underlying matrix is symmetric positive semi-definite (SPSD), the randomized \Ny{} decomposition is often preferred over the randomized SVD due to its lower computational cost, its single-pass nature and comparable approximation accuracy \cite{GittensMahoney2016, martinssonRandomizedNumericalLinear2020a, nakatsukasaFastStableRandomized2020}.
However, the restriction of the \Ny{} decomposition to SPSD matrices significantly limits its practical scope, as many matrices arising in data science and scientific computing are generally non-symmetric or indefinite, such as weight matrices in machine learning, non-self-adjoint discretized PDE operators, or Jacobians in nonlinear optimization. 
For matrices that are not SPSD, the \GN{} decomposition may be used instead \cite{nakatsukasaFastStableRandomized2020}, but
until now, there has not been a fast method to estimate the approximation error with a leave-one-out estimator.

In this work, we derive three fast randomized methods of leave-one-out error estimation for \GN{} approximations, corresponding to three different problem environments.
These formulas enable rank-adaptive algorithms that build accurate and cost-effective generalized \Ny{} approximations when the target rank for the approximation is unknown \emph{a priori}. 
We compare these three approaches in several numerical experiments with matrices designed to be adversarial and matrices that arise in practical applications, demonstrating that two of the methods (the leave-right-out and leave-twins-out (when applicable)) are most effective.
Our code is available at https://github.com/lorenzolazzarino/loo-estimators-for-gn.

 \section{Preliminaries}
 \label{sec:bg}
We summarize notation and key concepts used in our work. 

\subsection{Notation}
\label{sec:notation}

For a matrix $X \in \F^{m \times n}$, we  denote the $i,j$ entry of $X$ by $\jk{X}{i}{j}$.
We let $\jcol{X}{j}$ denote the $j^\text{th}$ column of $X$, and $\jrow{X}{i}$ denote the $i^\text{th}$ row. 
We will also frequently refer to matrices without certain rows or columns; to that end, let $\njcol{X}{j}$ be the matrix $X$ without the $j^\text{th}$ column,  $\njrow{X}{i}$ be the matrix $X$ without the $i^\text{th}$ row, and $\njk{X}{i}{j}$ be the matrix $X$ without the $i^\text{th}$ row and $j^\text{th}$ column. Given a matrix $X \in \F^{m \times n}$, we denote its $2$-norm by $\|X\|$, and its Frobenious norm by $\|X\|_F$.

We let $X^*$ denote the conjugate transpose of $X$, $X^\dagger$ its Moore-Penrose pseudoinverse, and we say a matrix $U$ is \textit{orthonormal} if its columns are orthonormal so that $U^* U = I$, the identity matrix.
The vector of all zeros is denoted by $\underline{0}$. 
The \Ny{} approximation of a matrix $X$ is denoted by $\widetilde X$, the \GNa{} by $\widehat{X}$, and the Schur complement corresponding to the $i,j$ entry of $X$ is denoted by $X / \jk{X}{i}{j} = \njk{X}{i}{j} - X_{-i,j} \jk{X}{i}{j}^\dagger X_{i,-j}$.

\subsection{Matrix decompositions}
\label{sec:matdecomps}

We briefly review two fundamental matrix decompositions frequently utilized in our work: the QR decomposition and the singular value decomposition (SVD).

\paragraph{The QR Decomposition}
Every matrix $X \in \F^{m \times n}$ admits a {QR factorization} $X = QR$,
where $Q \in \F^{m \times d}$ is orthonormal, $R \in \F^{d \times n}$ is upper triangular, $d = \min(m,n)$.
\vspace{2mm}

\paragraph{The Singular Value Decomposition}
Every matrix $X \in \F^{m \times n}$ admits a {singular value decomposition} (SVD), given by $X = U \Sigma V^*$, matrices $U$ and $V$ are orthonormal, and $\Sigma$ is diagonal.
The columns $\{ u_i \}_{i=1}^{d}$ and $\{v_{i}\}_{i=1}^{d}$ of $U$ and $V$ are called the left and right singular vectors of $X$, where $d = \min(m,n)$.
The diagonal elements $\{ \sigma_i\}_{i=1}^{d}$ of $\Sigma$ are the singular values of $X$, ordered so that $\sigma_1 \geq \sigma_2 \geq \cdots \geq \sigma_d \geq 0$.
A rank-$k$ truncated SVD is given by $X_k = \sum_{i=1}^k \sigma_i u_i v_i^*$.

\subsection{Randomized rangefinding}
\label{sec:randsketch}

For efficient low-rank approximations, particularly when only a matrix-vector primitive is available, randomized sketching  has proven to be a highly effective tool for computing fast and accurate low-rank approximations.
Randomized sketching allows us to approximate the column  space of a given matrix $X$ by analyzing how $X$ acts on random vectors, or more generally on matrices drawn from random distributions. 
For $X \in \F^{m \times n}$, the product $X\Omega \in \F^{m \times s}$ is known as a random (column) sketch of $X$, where $\Omega \in \F^{n \times s}$ is a random matrix with $s \ll n$ in general. 
Of note, the columns of $X\Omega$ are random linear combinations of the columns of $X$. 
A random row sketch $X^* \Phi$ (or $\Phi^*X$) for random matrix $\Phi \in \F^{m \times s}$ is analogous.

When $\Omega$ is a random matrix, such as one with entries drawn from a standard Gaussian distribution, the action of $\Omega$ preserves the geometry of the column space of $X$ with high probability; the column space of the sketch $X\Omega$ closely approximates that of $X$. 
This is critical for the performance of randomized algorithms that rely on the \emph{randomized rangefinder procedure}, cf. \cite[Algorithm 4.1]{halkoFindingStructureRandomness2011}, \cite[Algorithm 7]{martinssonRandomizedNumericalLinear2020a}.
The randomized rangefinder procedure computes an orthonormal matrix $Q$ such that $\| X - Q Q^* X \|$ is small, i.e. the columns of $Q$ form an approximate basis for the column space of $X$.
This is accomplished by the following steps:
\begin{enumerate}
    \item Draw $n \times s$ random matrix $\Omega$ with $s > k$,
    \item Form the $m \times s$ sketch $Y = X \Omega$,
    \item Compute $m \times k$ matrix $Q$ via rank-$k$ truncated QR decomposition $Y = QR$.
\end{enumerate}
In particular, if $\Omega$ is a Gaussian random matrix, a sketch size of $s = k + p$ with $p = O(1)$ is sufficient (e.g. $p=5$) for an accurate approximation with high probability, cf. \cite[Theorem 11.5]{martinssonRandomizedNumericalLinear2020a}.
This is often true in practice for other structured random matrices, e.g. subsampled randomized trigonometric transforms \cite{troppImprovedAnalysisSubsampled2011} or sparse sign matrices \cite{tropp2017practical}, though they are not as easily analyzable as Gaussian matrices.

\subsection{\Ny{} and \GNa s}
\label{sec:N-GN-approx}

The \Ny{} approximation is another low-rank approximation of an SPSD matrix $X \in F^{m \times n}$.
For arbitrary matrix $S \in \F^{n \times s}$, the \Ny{} approximation of $X$ with respect to $S$ is given by the SPSD matrix 
\begin{align}
    \label{eq:Nystrom}
    \widetilde{X} \langle S \rangle \equiv \widetilde{X} = \left (XS \right ) \left (S^* X S \right )^{-1} \left (XS \right )^*.
\end{align}
The error in the \Ny{} approximation is given by the Schur complement:
\begin{align}
    \label{eq:Nystrom_error}
    X - \widetilde{X}\langle S \rangle = X/S.
\end{align}
If $S$ is taken to be a Gaussian sketching matrix, then the rank-$k$ \Ny{} approximation $\widetilde{X}$ is nearly as accurate as the optimal rank-$k$ approximation, again with just slight oversampling, e.g., $s = k+5$.

The \GNa{} extends  (\ref{eq:Nystrom}) to non-SPSD matrices. 
Given any matrix $X \in \F^{m \times n}$, we now consider two arbitrary matrices $S \in \F^{n \times s}$ and $T \in \F^{m \times r}$.
The \GNa{} is given by
\begin{align}
    \label{eq:genNys}
    \widehat{X} \langle S,T\rangle \equiv \widehat{X} = \left (X S \right ) \left ( T^* X S \right )^\dag \left ( T^*X \right ).
\end{align}
If $S$ and $T$ are taken to be Gaussian sketching matrices, then the approximation error of $\widehat{X}$ is nearly optimal, cf. \cite[Theorem 3.2]{nakatsukasaFastStableRandomized2020}. The number of columns, $s$ and $r$, in the matrices $S$ and $T$ need not be equal. 
In fact, it has been suggested that introducing a discrepancy between $s$ and $r$ can enhance the accuracy of the \GNa~\cite{nakatsukasaFastStableRandomized2020}.
Throughout our work, we refer to the case where $r = s$ as the non-discrepant case, and the case where $r\geq s$ as the discrepant case.
For both \Ny{} and \GN, we denote the core matrix, $S^*XS$, or $\Omega^*X \Omega$ for random test matrix $\Omega$, and analogously $T^* X S$, or $\Phi^* X \Omega$, by $H$. Throughout the paper, we assume that $H$ is full-rank.


\subsection{Leave-one-out error estimation}
\label{sec:LOO-existing}
When computing a \Ny{} approximation, $\widetilde A$, of a given symmetric positive-definite matrix, $A$, using the random matrix $\Omega$, we first form the sketched matrix $A\Omega$ and then apply the inverse of the core matrix $H$ to one of the outer factors via its QR factorization. 
This procedure leverages the randomized rangefinder of Section~\ref{sec:randsketch}. 
We note that in general, it suffices to represent the \Ny{} approximation as in (\ref{eq:Nystrom}); however, there is an alternative form of the decomposition that can be computed via Cholesky decomposition of $H$ and SVD of the core matrix; see, e.g., \cite[Algorithm 4.1]{EpperlyTropp_VarianceEst2024}. 
Regardless, computing $\|A - \widetilde{A}\|$ na\"ively can undermine the computational benefits afforded by randomization.
In particular, in matrix-free or streaming environments, approximating the error using only the information needed to compute $\widetilde{A}$ is particularly desirable. 

We can efficiently estimate $\|A - \widetilde{A}\|$ by computing the \textit{leave-one-out (LOO) error}, a type of cross-validation. 
In cross-validation approaches, the data is partitioned, and, in each partition, some data is allocated into a training set and the rest into a testing set. 
For each partition, the model is fitted on the training set, and its error evaluated using the corresponding test set. 
The cross-validation error estimator is then computed by averaging the error for each partition \cite{wasserman2010alla}, providing
a nearly unbiased estimator of the true model error.
In particular, the estimator's bias decreases as the size of the training set approaches the number of data points \cite{celisse2008model}. 
The LOO approach is the most extreme case of cross-validation, where we allocate only one observation to the testing set and the rest to the training set. 
While the LOO approach typically provides the least-biased estimates  of $\|A - \widetilde{A}\|$, cf. \cite[Theorem 20.6]{wasserman2010alla}, it is also the most expensive, as, for a \Ny{} approximation with $\Omega \in \mathbb{F}^{n \times s}$, each error estimate requires an additional \Ny{} computation with $s-1$ vectors. 

A fast formula to compute the LOO estimator for \Ny{} is derived in \cite{EpperlyTropp_VarianceEst2024}, in which no additional randomized \Ny{} approximations are needed beyond $\widetilde{A}$, using random test matrix $\Omega$ drawn from an isotropic distribution, e.g., Gaussian. 
Let $\omega_j$ denote the $j$th column of $\Omega$, and let $\jcol{\Omega}{-j}$ be $\Omega$ with $\omega_j$ removed.
For $j=1,\ldots,s$, let ${\widetilde{A}}_{-j}$ be the rank-$(s-1)$ \Ny{} decomposition computed from $\jcol{\Omega}{-j}$, called a \textit{replicate}.
The {LOO error estimator} is
\begin{align}
\label{eq:naiveLOO}
    \textup{LOO} := \frac{1}{\sqrt{s}} \sum_{j=1}^s  \left \| \left ( A - {\widetilde{A}}_{-j} \right ) \omega_j \right \|_2^2,
\end{align}
which is an unbiased estimator for the mean-square error of the rank-$(s-1)$ \Ny{} approximation \cite[Theorem 2.1]{EpperlyTropp_VarianceEst2024}, i.e., $ \mathbb{E} \left \| A -  {\widetilde{A}}_{-s} \right \|^2_F = \textup{LOO}.$
To compute the replicates ${\widetilde{A}}_{-j}$ in (\ref{eq:naiveLOO}) efficiently, we use the following update formula from \cite[Eq. (4.1)]{EpperlyTropp_VarianceEst2024}, \cite[Eq. (2.4)]{EpperlyTroppWebber_XTrace2024} that depends only on pre-computed quantities:
\begin{align}
    \label{eq:LOOreps}
     {\widetilde{A}}_{-j} = A - \frac{t_j t_j^*}{\jk{H^{-1}}{j}{j}}, \qquad j = 1,\ldots,s,
\end{align}
where $t_j$ is the $j$-th column of $T = \begin{bmatrix}
    t_1 \cdots t_s
\end{bmatrix} =    A\Omega H ^{-1}$. 
Thus, given QR factorization $A\Omega = QR$, we can rewrite \eqref{eq:naiveLOO} using \eqref{eq:LOOreps}, obtaining the  replicate formula:
\begin{equation}
    \label{eq:replicateLOO}
    \textup{LOO} = \frac{1}{\sqrt{s}} \left\Vert R H^{-1} \textup{diag}\left\{ \frac{1}{\jk{H^{-1}}{j}{j}} \ : \ j=1,\dots,s \right\} \right\Vert_F
\end{equation}
\section{Fast estimation of approximation error for \GN}
\label{sec:Loo-GN}

Our main results are efficient error estimation techniques for the \GNa{} of a matrix $A$ using isotropic random matrices $\Omega$ and $\Phi$. 
We call the $j$-th column $\omega_j$ of $\Omega$ a right sample vector, and the $\ell$-th column $\phi_\ell$ of $\Phi$ a left sample vector. For $j=1,\dots ,s$ and $\ell = 1,\dots, r$, let $\njk{\widehat{A}}{\ell}{j}$ be the rank-$(s-1)$ \GNa{} computed by omitting those sample vectors.

Since in \GN{} we have two sample sets, there are multiple ways we can generalize the leave-one-out approach: the leave-pair-out (LPO), the leave-twins-out (LTO) and the leave-right-out (LRO) error estimators.
We first consider the non-discrepant case $r=s$, where replicate formulas for LPO and LTO are available, then derive our LRO formula in the discrepant case.
To illustrate the need for fast  replicate formulas, we first show how they can be na\"ively computed.

\subsection{The LPO and LTO formulas}
We focus here on the non-discrepant case, where $r=s$.
We begin by observing that we can estimate the error by computing all possible \GN{} estimates for which any left sample vector and any right sample vector are omitted, which we call the leave-pair-out (LPO) formula:
\begin{equation}
\label{eq:NaiveLPO}
	\text{LPO} := \frac{1}{s} \left( \sum_{j=1}^{s} \sum_{\ell =1}^{s} \left| \phi_\ell^* (A - \njk{\widehat A}{\ell}{j} ) \omega_j \right|^2 \right)^{\frac{1}{2}}.
\end{equation}
We can also leave out pairs of sample vectors with the same index, yielding
\begin{equation}
\label{eq:NaiveLTO}
	\text{LTO} := \frac{1}{\sqrt{s}} \left( \sum_{j=1}^{s}  \left| \phi_j^* (A - \njk{\widehat A}{j}{j} ) \omega_j \right|^2 \right)^{\frac{1}{2}},
\end{equation}
which we call the leave-twin-out (LTO) formula. 
In both cases, computing the na\"ive formula is impractical, as it requires computation of a new \GNa{} for each left-out pair or twin. 
In other words, evaluating the error for a single \GNa{} necessitates computing $s^2$ (for LPO) or $s$ (for LTO) additional approximations, making evaluating the error more expensive than computing the approximation. 
Therefore, we need a simple and effective way of computing the LPO and LTO formulas that avoids the need to recompute a \GNa{}. 
We derive these fast formulas next. \\

\paragraph{Leave-Pair-out (LPO) replicate formula}
\label{subsec:pairs}
We now derive a replicate formula for (\ref{eq:NaiveLPO}). 
Namely, we express \eqref{eq:NaiveLPO} in terms of pre-computed quantities, in this case entries of the (inverse of) core matrix $H$:
\begin{equation}
\label{eq:replicateLPO}
\begin{aligned}
	\text{LPO} :&= \frac{1}{s} \left( \sum_{j,\ell =1}^{s} \left|  \frac{1}{\jk{H^{-1}}{j}{\ell}}\right|^2 \right)^{\frac{1}{2}} \\[1ex]
    &= \frac{1}{s}\left\Vert \left\{ \frac{1}{\jk{H^{-1}}{j}{\ell}} \ : \ j=1,\dots s, \ \ell = 1,\dots s \right\}\right\Vert_F.
    \end{aligned}
\end{equation}
Note that the matrix on the right-hand side is simply the element-wise reciprocal of $H^{-1}$. Note, $H$ is assumed to be full rank and, thus, since $r=s$, it is invertible.

To derive \eqref{eq:replicateLPO}, we first express the \GNa{} $\hat A$ in terms of an approximation where a left and right sample vector are left out. 
Without loss of generality, we focus on how $\hat A$ relates to $\njk{\hat A}{s}{s}$, i.e., the \GN{} decomposition computed after removing the last sample vectors.
Specifically, we aim to characterize the relationship between the core matrices $H$ and $\njk{H}{s}{s}$.

To this end, we first write $H$ as a $2\times2$ block matrix, 
\begin{equation}
H =  \begin{bmatrix}
    \njk{H}{s}{s} & \njcol{\Phi}{s}^*A\omega_s \\ \phi_s^*A\njcol{\Omega}{s} & \phi_s^*A\omega_s
\end{bmatrix} =: \begin{bmatrix}
    \Bar{H} & \hb{1}{2} \\ \hb{2}{1} & \hb{3}{2}
\end{bmatrix} .
\end{equation}
We next apply the Banachiewicz inversion formula \cite{zhangSchurComplementIts2005}, so that
\begin{equation}
\label{eq:Inversion-pairs}
H^{-1} = \begin{bmatrix} \Bar{H}^{-1} + \Bar{H}^{-1} \hb{1}{2}(H/\Bar{H})^{-1} \hb{2}{1}\Bar{H}^{-1} & -\Bar{H}^{-1} \hb{1}{2} (H/\Bar{H})^{-1} \\
- (H/\Bar{H})^{-1}\hb{2}{1} \Bar{H}^{-1} & (H/\Bar{H})^{-1} \end{bmatrix}
\end{equation}
where  $H/\Bar{H} = \hb{3}{2} - \hb{2}{1} \Bar{H}^{-1}\hb{1}{2} = \phi_s^*A\omega_s - (\phi_s^*A\njcol{\Omega}{s}) (\njk{H}{s}{s})^{-1} (\njcol{\Phi}{s}^*A\omega_s) \in \mathbb{R}$ is the Schur complement. 
We then have
\begin{equation}
	\label{eq:finalInversion-pairs}
	H^{-1} - \frac{H^{-1} e_s e_s^* H^{-1}}{e_s^* H^{-1} e_s} = \begin{bmatrix} \njk{H}{s}{s}^{-1} & \underline{0} \\ \underline{0}^* & 0 \end{bmatrix},
\end{equation}
i.e., the desired relation that allows us to express $\njk{\widehat A}{s}{s}$ in terms of $\widehat A$ as follows:
\begin{equation}
	\label{eq:firstXrs}
	\begin{aligned}
		\widehat A &= A\Omega H^{-1} \Phi^*A \\
		  &= A\Omega \left( \begin{bmatrix} \njk{H}{s}{s}^{-1} & \underline{0} \\ \underline{0}^* & 0 \end{bmatrix} + \frac{H^{-1} e_s e_s^* H^{-1}}{e_s^* H^{-1} e_s} \right) \Phi^*A \\
		  &= A\njcol{\Omega}{s} \njk{H}{s}{s}^{-1} (\njcol{\Phi}{s})^*A + \frac{A\Omega H^{-1} e_s e_s^* H^{-1} \Phi^*A}{e_s^* H^{-1} e_s} \\
          &= \njk{\widehat A}{s}{s} + \frac{A\Omega H^{-1} e_s e_s^* H^{-1} \Phi^*A}{e_s^* H^{-1} e_s}.
	\end{aligned}
\end{equation}
Letting $T = \begin{bmatrix} t_1 \cdots t_s \end{bmatrix} = A\Omega H^{-1}$ and $Z = \begin{bmatrix} z_1 \cdots z_s \end{bmatrix} = A^*\Phi (H^{-1})^*$, we can write
\begin{equation}
	\label{eq:Xrs-pairs}
		\njk{\widehat A}{s}{s}  
		= \widehat A - \frac{t_s z_s^*}{\jk{H^{-1}}{s}{s}}.
\end{equation}
Before we use these equations to obtain a replicate formula, we quickly analyze
\begin{equation} \label{eq:replicateFirst}
	\begin{aligned}
		\left| \phi_\ell^* (A - \njk{\widehat A}{\ell}{j} ) \omega_j \right|^2 & = \left| \phi_\ell ^* (A - \widehat A + \frac{t_\ell z_j^*}{\jk{H^{-1}}{j}{\ell}}  ) \omega_j \right|^2 \\
									&= \left| \phi_\ell^* (I - A \Omega (\Phi^* A \Omega)^{-1} \Phi^*)A\omega_j + \phi_\ell^*\frac{t_\ell z_j^*}{\jk{H^{-1}}{j}{\ell}} \omega_j \right|^2.
		\end{aligned}
		\end{equation}
Since $A \Omega (\Phi^* A \Omega)^{-1} \Phi^*$ is an oblique projection\footnote{Note $(A \Omega (\Phi^* A \Omega)^{-1} \Phi^*)^2 = A \Omega (\Phi^* A \Omega)^{-1} \Phi^* A \Omega (\Phi^* A \Omega)^{-1} \Phi^* = A \Omega (\Phi^* A \Omega)^{-1} \Phi^*$. 
Also, if  $x = z + y$, $y \in \textbf{Range}(\Phi^*)$ and $z \in \textbf{Null}(\Phi)$, then $ A \Omega (\Phi^* A \Omega)^{-1} \Phi^* x =  A \Omega (\Phi^* A \Omega)^{-1} \Phi^* y \in \textbf{Range}(A\Omega)$. } 
onto the range of $A \Omega $, we have $(A \Omega (\Phi^* A \Omega)^{-1} \Phi^*)A\omega_j = A\omega_j$. 
Then \eqref{eq:replicateFirst} becomes
\begin{equation}
	\label{eq:replicateSecond}
\left| \phi_\ell^* (A - \njk{\widehat A}{\ell}{j} ) \omega_j \right|^2  
= \left| \frac{1}{\jk{H^{-1}}{j}{\ell}}  \phi_\ell^*A\Omega H^{-1} e_\ell e_j^* H^{-1} \Phi^*A \omega_j\right|^2.
\end{equation}
where  $\frac{1}{\jk{H^{-1}}{j}{\ell}}$ is a scalar.
Finally, we notice that
\begin{equation}
		 \phi_\ell^* A\Omega = e_\ell^* H, \quad
		  \Phi^*A\omega_j = He_j,
\end{equation}
and therefore, by~\eqref{eq:replicateSecond},
\begin{equation}
	\begin{aligned}
	\left| \phi_\ell^* (A - \njk{\widehat A}{\ell}{j} ) \omega_j \right|^2 & = \left| \frac{1}{\jk{H^{-1}}{j}{\ell}} e_\ell^* H H^{-1} e_\ell e_j^* H^{-1} H e_j \right|^2 
								 = \left| \frac{1}{\jk{H^{-1}}{j}{\ell}}\right|^2,
\end{aligned}
\end{equation}
from which the desired LPO replicate formula~\eqref{eq:replicateLPO} follows. \\

\paragraph{Leave-Twins-out (LTO) replicate formula}
\label{subsec:twins}
Considering all possible pairs, as in the LPO formula, may account for more summands in the error estimate than necessary in practice. Thus, we consider an analogous formula, where fewer index couples are taken into account.
Instead of looking at all possible pairs, as for the LPO formula, we can look only at pairs with the same index, or a twin index, as in the LTO formula~\eqref{eq:NaiveLTO}.
For a fast LTO formula, we consider  those terms in the sums in the LPO formula for which $j=\ell$, that is 
\begin{equation}
\label{eq:replicateLTO}
\text{LTO} := \frac{1}{\sqrt{s}} \left( \sum_{j=1}^{s} \left|  \frac{1}{\jk{H^{-1}}{j}{j}}\right|^2 \right)^{\frac{1}{2}} = \frac{1}{\sqrt{s}} \left\Vert \text{diag}\left\{ \frac{1}{\jk{H^{-1}}{i }{i}} : i = 1,\dots,s\right\} \right\Vert_F
\end{equation}

\subsection{The LRO formula}
 For the more general case when $r \geq s$, the core matrix $H$ is rectangular with full column rank, and the Banachiewicz inversion formula no longer holds. 
 This results in the replicate formulas of the previous section no longer holding. To remedy this, we now introduce the leave-right-out (LRO) formula. 
 Here, the error is estimated by leaving out a right sample vector, that is:
\begin{equation}
\label{eq:NaiveLRO}
	\text{LRO} := \frac{1}{\sqrt{s}} \left( \sum_{j=1}^{s}  \left\Vert (A - \njcol{\widehat A}{j} ) \omega_j \right\Vert^2 \right)^{\frac{1}{2}}.
\end{equation}
As before, the LRO formula as in \eqref{eq:NaiveLRO} is not practical, requiring $s$ \GNa s, leading us to derive the following fast formula for \eqref{eq:NaiveLRO}:
\begin{equation}
\label{eq:replicateLRO}
\begin{aligned}
	\text{LRO} &= \frac{1}{\sqrt{s}} \left( \sum_{j=1}^{s}  \left\Vert R_\Omega \frac{\jcol{(H^*H)^{-1}}{j}}{\jk{(H^*H)^{-1}}{j}{j}} \right\Vert^2 \right)^{\frac{1}{2}}\\[1ex]
    &=\frac{1}{\sqrt{s}}\left\Vert R_\Omega (H^*H)^{-1} \text{diag}\left\{ \frac{1}{\jk{(H^*H)^{-1}}{i }{i}} : i = 1,\dots,s\right\} \right\Vert_F,
    \end{aligned}
\end{equation}
where, $R_\Omega$ is the $R$-factor of the $QR$-factorization of $A\Omega$, i.e., $A\Omega = Q_\Omega R_\Omega$.

To derive this fast replicate formula, we again characterize the relationship between (the pseudoinverses of) $H$ and $\njcol{H}{s}$. 
While the Banachiewicz inversion formula \eqref{eq:Inversion-pairs} does not hold for pseudoinverses, we use the observation from \cite{guttel2024shermanmorrisonwoodburyapproachsolving} that, since $H$ (and $\njcol{H}{s}$) is full column rank, it is possible to write its pseudoinverse 
as\footnote{The case where $s\geq r$ is analogous: $H$ would be full row rank, and $H^\dagger = H(HH^*)^{-1}$.}
:
\begin{equation}
    H^\dagger = (H^*H)^{-1}H^*.
\end{equation}
We can then obtain the equivalent of equation \eqref{eq:finalInversion-pairs} for $H^*H$, that is,
\begin{equation}
    (H^*H)^{-1} - \frac{(H^*H)^{-1}e_s e_s^*(H^*H)^{-1}}{e_s(H^*H)^{-1}e_s}= \begin{bmatrix}
     (\njcol{H}{s}^*\njcol{H}{s})^{-1}   & 0 \\ 0 & 0
    \end{bmatrix}.
\end{equation}
Multiplying on the right by $H^*$, we obtain
\begin{equation}
    H^\dagger - \frac{(H^*H)^{-1}e_s e_s^*H^\dagger}{e_s(H^*H)^{-1}e_s}= \begin{bmatrix}
     \njcol{H}{s}^\dagger   \\ 0 
    \end{bmatrix}.
\end{equation}
Following similar steps as for LPO, we can relate $\widehat A$ to $\njcol{\widehat A}{s}$:
\begin{equation}
\begin{aligned}
    \widehat A &= A\Omega H^\dagger \Phi^*A \\
    &= A\Omega \left(  \begin{bmatrix}
     \njcol{H}{s}^\dagger   \\ 0 
    \end{bmatrix}+\frac{(H^*H)^{-1}e_s e_s^*H^\dagger}{e_s(H^*H)^{-1}e_s}\right)\Phi^*A \\
    &= A\njcol{\Omega}{s} \njcol{H}{s}^\dagger \Phi^*A + \frac{A\Omega(H^*H)^{-1}e_s e_s^*H^\dagger\Phi^*A}{e_s(H^*H)^{-1}e_s}.
\end{aligned}
\end{equation}
That is, we can write
\begin{equation}
    \njcol{\widehat A}{s} = \widehat A - \frac{f_s g_s}{\jk{(H^*H)^{-1}}{s}{s}}
\end{equation}
where $f_s$ and $g_s$ are the last column of the matrices $F = \begin{bmatrix} f_1 \cdots f_s \end{bmatrix} = A\Omega (H^*H)^{-1}$ and $G = \begin{bmatrix} g_1 \cdots g_s \end{bmatrix} = A^*\Phi (H^\dagger)^*$.
We observe that
\begin{equation}
\label{eq:replicateFirstLRO}
    \begin{aligned}
        \left\Vert (A- \njcol{\widehat A}{j})\omega_j \right\Vert^2 &= \left\Vert (A- \widehat A - \frac{f_j g_j}{\jk{(H^*H)^{-1}}{j}{j}})\omega_j \right\Vert^2 \\
        &=\left\Vert A\omega_j- \widehat A \omega_j- \frac{f_s g_s}{\jk{(H^*H)^{-1}}{j}{j}}\omega_j \right\Vert^2 \\
        &=\left\Vert (I - A\Omega(\Phi^*A\Omega)^\dagger\Phi^*)A \omega_j- \frac{f_j g_j}{\jk{(H^*H)^{-1}}{j}{j}}\omega_j \right\Vert^2. \\
    \end{aligned}
\end{equation}
Since $A \Omega (\Phi^* A \Omega)^{-1} \Phi^*$ is an oblique projection onto the column space of $A \Omega $, we have $(A \Omega (\Phi^* A \Omega)^{-1} \Phi^*)A\omega_j = A\omega_j$, and~\eqref{eq:replicateFirstLRO} becomes
\begin{equation}
\label{eq:single-termLRO}
    \begin{aligned}
        \left\Vert (A- \njcol{\widehat A}{j})\omega_j \right\Vert^2 &=
        \left\Vert  \frac{f_j g_j}{\jk{(H^*H)^{-1}}{j}{j}}\omega_j \right\Vert^2 \\
        &= \left\Vert \frac{1}{\jk{(H^*H)^{-1}}{j}{j}} A\Omega(H^*H)^{-1}e_je_j^*H^\dagger \Phi^*A\omega_j \right\Vert^2 \\
        &= \left\Vert \frac{1}{\jk{(H^*H)^{-1}}{j}{j}} A\Omega(H^*H)^{-1}e_je_j^*H^\dagger H e_j \right\Vert^2 \\
        &= \left\Vert \frac{A\Omega (H^*H)^{-1}e_j}{\jk{(H^*H)^{-1}}{j}{j}}\right\Vert^2,
    \end{aligned}
\end{equation}
where the last equality follows from $H$ having full column rank, so $H^\dagger H = I$. Then, by considering the $QR$-factorization, $A\Omega = Q_\Omega R_\Omega$, we can use the unitary invariance of the $2$-norm, and \eqref{eq:single-termLRO} in \eqref{eq:NaiveLRO}, to obtain the desired replicate formula \eqref{eq:replicateLRO}.

The na\"ive and the fast replicate formulas for the leave-right-out estimator hold for $r\geq s$. In the next section, we numerically compare all of our fast formulas, against their na\"ive formulas, as well as the true \GN{} error.
\section{Numerical Experiments}

\label{sec:NumExp} 
In this section, we analyze the behavior of the presented estimators through numerical experiments. For each experiment, we compute the true error of the \GNa{} for increasing target ranks~$s$, compared to the estimates provided by our leave-pair-out (LPO), leave-twin-out (LTO), and leave-right-out (LRO) estimators. 
Through three experiments, we compare the runtimes and accuracies of our na\"ive and fast replicate formulations. 
We first show a favorable example of a matrix with an exponentially-decaying spectrum, then move on to an adversarial matrix.
We also show a practical example from advection-diffusion equations. 
In each experiment \footnote{A detailed implementation can be found at https://github.com/lorenzolazzarino/loo-estimators-for-gn.}, we consider matrices of size $500\times500$, and we compute the \GNa{} for increasing target rank~$s$, starting from $s=25$ and incrementing by $25$ up to $s=250$.

\vspace{2mm}
\paragraph{{Exponentially decaying singular values}}
We use a classical example in which the \GNa{} performs well, namely, when the approximated matrix has exponentially decaying singular values, to compare the na\"ive and fast formulations of each proposed estimator.
We consider a matrix whose singular values are defined as $\sigma_i = 2^{-\frac{i}{6}}$. 
The matrices of singular vectors are generated as Haar matrices, i.e., the orthogonal factors of Gaussian matrices.  We perform experiments in both the discrepant ($r= s +5$) and non-discrepant ($r=s$) cases. 

Figures~\ref{fig:exp} and \ref{fig:exp-discrepant} show the results of a single run; in all three cases, numerical stability does not affect the replicate formulas, which remain numerically equivalent to their na\"ive counterparts. 
We observe rapid growth in the runtime for the na\"ive formulas, especially for LPO, which requires substantially more computations than the other formulas.
To emphasize this, we report the computational runtime as well.
\begin{figure}
    \centering
    \begin{tikzpicture}[scale = \boxplotscale, anchor = west]
	\pgfplotstableread[col sep=comma]{./csvs/output_brute_exp6_m500_n500_os0.csv}\data
    \pgfplotstableread[col sep=comma]{./csvs/output_brute_exp6_m500_n500_os5.csv}\datad
    \newcommand{\Matname}{Exp Matrix }
    \pgfplotscreateplotcyclelist{cross_colors}{cGN, cLTO, cLTO, cLPO, cLPO, cLCO, cLCO, cSVD}
    \begin{axis}[
        name = axis1,
        ylabel = {Froebinius Norm Error},
        xlabel = {Rank of Approximation},
        title = {Error},
        ymode = log,
        ymajorgrids,
        scale only axis,
        xlabel near ticks,
        ylabel near ticks,
        legend columns = -1,
        legend to name = lgd:exp6os0,
        legend style ={font = \tiny},
        ]        
        \addplot[gn_line]table[x = ranks, y = err_GN]\data;
        \addplot[lto_line]table[x = ranks, y = LTO]\data;
        \addplot[blto_line]table[x = ranks, y = bruteLTO]\data;
        \addplot[lpo_line]table[x = ranks, y = LPO]\data;
        \addplot[blpo_line]table[x = ranks, y = bruteLPO]\data;
        \addplot[lco_line]table[x = ranks, y = LCO]\data;
        \addplot[blco_line]table[x = ranks, y = bruteLCO]\data;
        \addplot[svds_line]table[x = ranks, y = err_svd]\data;
        \legend{GN, LTO, BruteLTO, LPO, BruteLPO, LRO, BruteLRO, SVD};
    \end{axis}
    \path (axis1.outer north east) ++(0.5, 0)  coordinate (bb);
    \begin{axis}[
        name = axis2,
        at={(bb)},
        anchor=outer north west,
       ylabel = {Froebinius Norm Error},
        xlabel = {Rank of Approximation},
        title = {Time},
        ymode = log,
        ymajorgrids,
        scale only axis,
        xlabel near ticks,
        ylabel near ticks,
        ]
        \addplot[gn_line]table[x = ranks, y = timeGN]\data;
        \addplot[lto_line]table[x = ranks, y = tLTO]\data;
        \addplot[blto_line]table[x = ranks, y = tbruteLTO]\data;
        \addplot[lpo_line]table[x = ranks, y = tLPO]\data;
        \addplot[blpo_line]table[x = ranks, y = tbruteLPO]\data;
        \addplot[lco_line]table[x = ranks, y = tLCO]\data;
        \addplot[blco_line]table[x = ranks, y = tbruteLCO]\data;
    \end{axis}
    \path (axis1.outer south west) ++(-3, -0.5)  coordinate (lgp);
    \node at (lgp) {\pgfplotslegendfromname{lgd:exp6os0}};
\end{tikzpicture}
     \caption{Non-discrepant case for exponentially decaying singular values. (Left) Error from the \GNa{}, with the na\"ive (brute) computations (solid) lines and fast replicate computations (dashed). (Right) Runtime plot.}
    \label{fig:exp}
\end{figure}
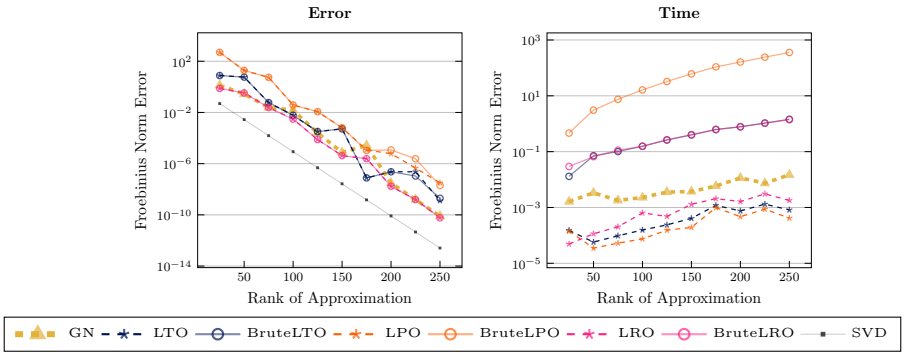
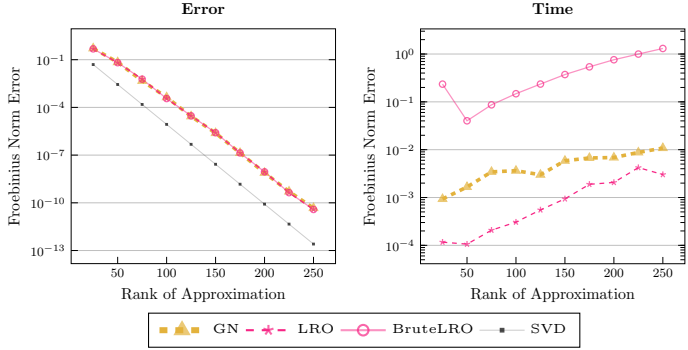
\begin{figure}
    \centering
    \begin{tikzpicture}[scale = \boxplotscale, anchor = west]
	\pgfplotstableread[col sep=comma]{./csvs/output_brute_exp6_m500_n500_os0.csv}\data
    \pgfplotstableread[col sep=comma]{./csvs/output_brute_exp6_m500_n500_os5.csv}\datad
    \newcommand{\Matname}{Exp Matrix }
    \pgfplotscreateplotcyclelist{cross_colors}{cGN, cLTO, cLTO, cLPO, cLPO, cLCO, cLCO, cSVD}
    \begin{axis}[
        name = axis1,
        ylabel = {Froebinius Norm Error},
        xlabel = {Rank of Approximation},
        title = {Error},
        ymode = log,
        ymajorgrids,
        scale only axis,
        xlabel near ticks,
        ylabel near ticks,
        legend columns = -1,
        legend to name = lgd,
        legend style ={font = \tiny},
        ]        
        \addplot[gn_line]table[x = ranks, y = err_GN]\datad;
        \addplot[lco_line]table[x = ranks, y = LCO]\datad;
        \addplot[blco_line]table[x = ranks, y = bruteLCO]\datad;
        \addplot[svds_line]table[x = ranks, y = err_svd]\datad;
        \legend{GN, LRO, BruteLRO, SVD};
    \end{axis}
    \path (axis1.outer north east) ++(.5, 0)  coordinate (bb);
    \begin{axis}[
        name = axis2,
        at={(bb)},
        anchor=outer north west,
       ylabel = {Froebinius Norm Error},
        xlabel = {Rank of Approximation},
        title = {Time},
        ymode = log,
        ymajorgrids,
        scale only axis,
        xlabel near ticks,
        ylabel near ticks,
        ]
        \addplot[gn_line]table[x = ranks, y = timeGN]\datad;
        \addplot[lco_line]table[x = ranks, y = tLCO]\datad;
        \addplot[blco_line]table[x = ranks, y = tbruteLCO]\datad;
    \end{axis}
    \path (axis1.outer south west) ++(3, -0.5)  coordinate (lgp);
    \node at (lgp) {\pgfplotslegendfromname{lgd}};
\end{tikzpicture}
     \caption{Discrepant case for exponentially decaying singular values. (Left) Error from the \GNa{} with the na\"ive (brute) computations (solid) lines and fast replicate computations (dashed). (Right) Runtime plot.}
    \label{fig:exp-discrepant}
\end{figure}

\vspace{2mm}
\paragraph{{Chan matrix}} For the next experiment, we consider the Chan matrix \cite{CHAN198767}, which has entries equal to 1 on the diagonal, $-1$ above, and 0 elsewhere; this is highly adversarial for low-rank approximation algorithms.
We again consider both the discrepant ($r= s +5$) and non-discrepant ($r=s$) cases, and compute both replicate and na\"ive, or brute, formulas for LTO and LPO (when applicable), as well as for LRO. 
We also plot the optimal error achievable for a rank-$s$ approximation, i.e., the error of rank-$s$ truncated SVD, to highlight the badly-behaved spectrum of the Chan matrix.

\begin{figure}
    \centering
    \begin{tikzpicture}[scale = \boxplotscale, anchor = west]
	\pgfplotstableread[col sep=comma]{./csvs/output_brute_chan_m500_n500_os0.csv}\data
    \pgfplotstableread[col sep=comma]{./csvs/output_brute_chan_m500_n500_os5-2.csv}\datad
    \newcommand{\Matname}{Chan Matrix }
    \pgfplotscreateplotcyclelist{cross_colors}{cGN, cLTO, cLTO, cLPO, cLPO, cLCO, cLCO, cSVD}
    \begin{axis}[
        name = axis1,
        ylabel = {Froebinius Norm Error},
        xlabel = {Rank of Approximation},
        title = {Error (Non-Discrepant)},
        ymode = log,
        ymajorgrids,
        scale only axis,
        xlabel near ticks,
        ylabel near ticks,
        legend columns = -1,
        legend to name = lgd:chanos0,
        legend style ={font = \tiny},
        ]        
        \addplot[gn_line]table[x = ranks, y = err_GN]\data;
        \addplot[lto_line]table[x = ranks, y = LTO]\data;
        \addplot[blto_line]table[x = ranks, y = bruteLTO]\data;
        \addplot[lpo_line]table[x = ranks, y = LPO]\data;
        \addplot[blpo_line]table[x = ranks, y = bruteLPO]\data;
        \addplot[lco_line]table[x = ranks, y = LCO]\data;
        \addplot[blco_line]table[x = ranks, y = bruteLCO]\data;
        \addplot[svds_line]table[x = ranks, y = err_svd]\data;
        \legend{GN, LTO, BruteLTO, LPO, BruteLPO, LRO, BruteLRO, SVD};
    \end{axis}
    \path (axis1.outer north east) ++(0.5, 0)  coordinate (bb);
    \begin{axis}[
        name = axis2,
        at={(bb)},
        anchor=outer north west,
       ylabel = {Froebinius Norm Error},
        xlabel = {Rank of Approximation},
        title = {Error (Discrepant)},
        ymode = log,
        ymajorgrids,
        scale only axis,
        xlabel near ticks,
        ylabel near ticks,
        ]
       \addplot[gn_line]table[x = ranks, y = err_GN]\datad;
        \addplot[lco_line]table[x = ranks, y = LCO]\datad;
        \addplot[blco_line]table[x = ranks, y = bruteLCO]\datad;
        \addplot[svds_line]table[x = ranks, y = err_svd]\datad;
     
    \end{axis}
    \path (axis1.outer south west) ++(-3, -0.5)  coordinate (lgp);
    \node at (lgp) {\pgfplotslegendfromname{lgd:chanos0}};
\end{tikzpicture}
    \caption{Error for the non-discrepant case (left)  and discrepant case (right) of the \GNa{} of the Chan matrix.}
    \label{fig:chan}
\end{figure}
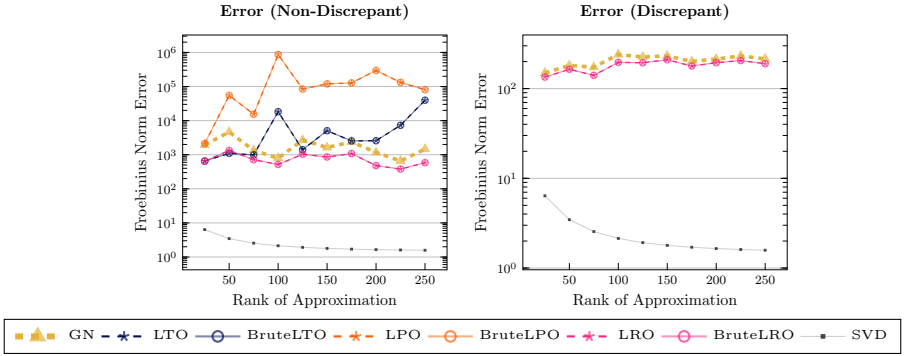

Figure~\ref{fig:chan} shows that the LTO (when applicable) and LRO formulas are very accurate in predicting the true \GN{} error, whereas the LPO formula can be less accurate. 
This is a desirable outcome, as the experiments indicate that the less expensive estimators provide better accuracy.

\vspace{2mm}
\paragraph{{Advection-diffusion}}
In practice, the \GN{} decomposition may be used for the approximation of Green's functions when an operator is not self-adjoint, 
as in the advection-diffusion equations:
\begin{equation}
\label{eq:advec}
    \frac{\partial C}{\partial t} =  \nabla \cdot (-\alpha \nabla C + u C) + f,
\end{equation}
where $\nabla$ is the gradient, $\nabla \cdot$ is the divergence operator, $\alpha$ is the diffusivity, $u$ is the velocity field, and $f$ is the source. 

\begin{figure}
    \centering
    \begin{tikzpicture}[scale = \boxplotscale, anchor = west]
	\pgfplotstableread[col sep=comma]{./csvs/addvec_diff.csv}\data
    \pgfplotstableread[col sep=comma]{./csvs/output_brute_exp6_m500_n500_os5.csv}\datad
    \newcommand{\Matname}{Exp Matrix }
    \pgfplotscreateplotcyclelist{cross_colors}{cGN, cLTO, cLTO, cLPO, cLPO, cLCO, cLCO, cSVD}
    \begin{axis}[
        name = axis1,
        ylabel = {Froebinius Norm Error},
        xlabel = {Rank of Approximation},
        title = {Error},
        ymode = log,
        ymajorgrids,
        scale only axis,
        xlabel near ticks,
        ylabel near ticks,
        legend columns = -1,
        legend to name = lgd:advec,
        legend style ={font = \tiny},
        ]        
        \addplot[gn_line]table[x = rank, y = GN_err]\data;
        \addplot[lto_line]table[x = rank, y = LTO]\data;
        \addplot[blto_line]table[x = rank, y = bruteLTO]\data;
        \addplot[lpo_line]table[x = rank, y = LPO]\data;
       \addplot[blpo_line]table[x = rank, y = bruteLPO]\data;
        \addplot[lco_line]table[x = rank, y = LCO]\data;
        \addplot[blco_line]table[x = rank, y = bruteLCO]\data;
        \addplot[svds_line]table[x = rank, y = svd_err]\data;
        \legend{GN, LTO, BruteLTO, LPO, BruteLPO, LRO, BruteLRO, SVD};
    \end{axis}
    \path (axis1.outer north east) ++(0.5, 0)  coordinate (bb);
    \begin{axis}[
        name = axis2,
        at={(bb)},
        anchor=outer north west,
       ylabel = {Froebinius Norm Error},
        xlabel = {Rank of Approximation},
        title = {Time},
        ymode = log,
        ymajorgrids,
        scale only axis,
        xlabel near ticks,
        ylabel near ticks,
        ]
        \addplot[gn_line]table[x = rank, y = tGN]\data;
        \addplot[lto_line]table[x = rank, y = tLTO]\data;
        \addplot[blto_line]table[x = rank, y = tbruteLTO]\data;
        \addplot[lpo_line]table[x = rank, y = tLPO]\data;
        \addplot[blpo_line]table[x = rank, y = tbruteLPO]\data;
        \addplot[lco_line]table[x = rank, y = tLCO]\data;
        \addplot[blco_line]table[x = rank, y = tbruteLCO]\data;
    \end{axis}
    \path (axis1.outer south west) ++(-3, -0.5)  coordinate (lgp);
    \node at (lgp) {\pgfplotslegendfromname{lgd:advec}};
\end{tikzpicture}
     \caption{Non-discrepant case for the advection-diffusion equations. (Left) Error from the \GNa{}, with the na\"ive (brute) computations (solid) lines and fast replicate computations (dashed). (Right) Runtime plot.}
    \label{fig:adv-dif}
\end{figure}
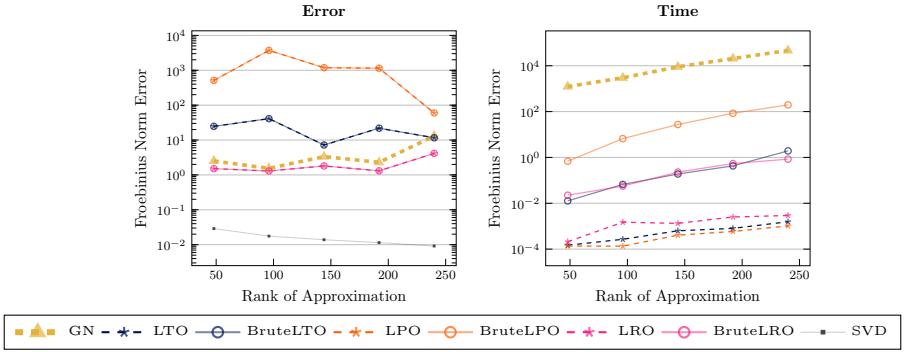

In our experiment, we solve the stationary problem associated with both \eqref{eq:advec} and its adjoint,
\begin{equation}
    \frac{\partial C}{\partial t}^* =  \nabla \cdot (-\alpha \nabla C - u C) + f,
\end{equation}
in two dimensions with Dirchlet boundary conditions. We discretize the system over an unit square with an 18 by 18 grid, and form a $324 \times 324$ \GNa{} to the Green's function associated with the steady state solution of (\ref{eq:advec}) with $\alpha = 1$ and $u = 1.5$. 
To ``sketch'' the Green's functions, we use the approach of \cite{boulle2022a, boulle2023learning}, generating random functions $f$ by drawing from a Gaussian Process distribution defined by the identity kernel \cite{boulle2024operator}. 
We again observe in Figure~\ref{fig:adv-dif} that the LRO estimator is effective and offers substantial speed-ups.
\section{Conclusion}
In this paper, we introduce three new leave-one-out error estimators for \GN{}. 
We derive computationally efficient replicate formulas that are cheaper alternatives to their na\"ive counterparts for the leave-pair-out (LPO), leave-twin-out (LTO), and leave-right-out (LRO) estimators. 
While the LPO and LTO estimators are applicable only in the absence of discrepancy, the leave-right-out (LRO) estimator applies to both the discrepant and non-discrepant settings.
Our numerical experiments underscore the effectiveness of these estimators. 
Notably, the LRO estimator, which leaves out the least amount of information, generally yields the most accurate estimates. 
However, when applicable, the leave-twin-out (LTO) estimator proves equally competitive.

This work serves as a foundational step toward developing error estimators for more complex (randomized) low-rank approximations. These concepts may be applicable in other problem frameworks, such as non-negative matrix factorizations, CUR factorizations, and a wide range of other (randomized) low-rank approximations.

\section*{Acknowledgments}
The authors would like to thank H. Al Daas for the fruitful discussions and invaluable suggestions, while writing this paper.

\bibliography{main_r}
\bibliographystyle{siamplain}
\end{document}